\documentclass[11pt]{amsart}

\usepackage{amssymb}
\usepackage{amsmath}
\usepackage{latexsym}
\usepackage{amsfonts}
\usepackage{graphicx}
\newtheorem{theorem}{Theorem}[section]

\newtheorem{corollary}[theorem]{Corollary}

\newtheorem{def-thm}[theorem]{Definition-Theorem}

\theoremstyle{definition}
\newtheorem{definition}[theorem]{Definition}

\newtheorem{remark}[theorem]{Remark}

\newtheorem{conjecture}[theorem]{Conjecture}

\allowdisplaybreaks

\newcommand{\be}{\begin{equation}}
\newcommand{\ee}{\end{equation}}
\newcommand{\bea}{\begin{eqnarray}}
\newcommand{\eea}{\end{eqnarray}}

\def\XXint#1#2#3{{\setbox0=\hbox{$#1{#2#3}{\int}$ }
\vcenter{\hbox{$#2#3$ }}\kern-.6\wd0}}

\begin{document}

\title{Positive curvature operator, projective manifold and rational connectedness}
\author{Kai Tang}
\address{Kai Tang, College of Mathematics and Computer Science, Zhejiang Normal University, Jinhua, 321004, Zhejiang, P.R. China}
\email{tangkai0810120925@163.com}

\maketitle

\markleft{Positive curvature operator and projective manifold} \markright{Positive curvature operator and projective manifold}

\begin{abstract}
In his recent work \cite{Y1}, X. Yang  proved a conjecture raised by Yau in 1982 (\cite{Yau82}), which states that any compact K\"{a}hler manifold with
positive holomorphic sectional curvature must be projective. In this note, we prove that any compact Hermitian manifold $X$ with positive real bisectional curvature, its hodge number $h^{1,0}=h^{2,0}=h^{n-1,0}=h^{n,0}=0$. In particular, if in addition $X$ is K\"{a}hler, then $X$ is projective. Also, it is rationally connected manifold when $n=3$. This partially confirms the conjecture 1.11 \cite{Y1} which is proposed by X. Yang.
\end{abstract}
\maketitle

\tableofcontents
\section{Introduction}
For compact K\"{a}hler manifold with positive holomorphic sectional curvature, S.-T. Yau proposed the following well-known conjecture \cite{Yau82}:
\vspace{0.2cm}

\noindent {\bf Conjecture}\,(Yau). {\em Let X be a compact K\"{a}hler manifold. If X has a
K\"{a}hler metric with positive holomorphic sectional curvature, then X is a
projective and rationally connected manifold.
}

\vspace{0.2cm}
A projective manifold $X$ is called rationally connected if any two points of $X$ can be connected by some rational curve. In his recent work \cite{Y1}, X. Yang confirmed the Yau's conjecture. More generally, he obtained

\begin{theorem} [{\cite{Y1}, Theorem1.7}] Let $(X,\omega)$ be a compact K\"{a}hler manifold with positive holomorphic sectional curvature, then for every $1\leq p\leq dimX$, $(\Lambda^{p}T_{X},\Lambda^{p}\omega)$ is RC-positive and $H^{p,0}_{\overline{\partial}}(X)=0$. In particular, $X$ is projective and rationally connected manifold.
\end{theorem}

It is worth noting that Heier-Wong also confirmed Yau's conjecture in the special case when X is projective \cite{HW}. The method of X. Yang crucially relied on the geometric properties of RC-positivity which he proved by using techniques in non-K\"{a}hler geometry, and also a
minimum principle for K\"{a}hler metrics with positive holomorphic sectional curvature, while Heier-Wong's method built on an
average argument and certain integration by parts on projective manifolds.

For Hermitian case, X. Yang also obtained
\begin{corollary}[{\cite{Y1}, Corollary 1.10}] Let $X$ be a compact K\"{a}hler surface. If there exists a Hermitian
metric $¦Ø$ with positive holomorphic sectional curvature, then $X$ is a
projective and rationally connected manifold.
\end{corollary}

When the complex dimension of manifold $n>2$, he proposed the following conjecture:

\begin{conjecture}[{\cite{Y1}, Conjecture 1.11}] {\em Let $X$ be a compact complex manifold of complex dimension $n>2$. Suppose $X$ has a Hermitian metric with positive holomorphic sectional
curvature, then $\Lambda^{p}T_{X}$ admits a smooth RC-positive metric for every $1<p<n$. In particular, if in addition $X$ is K\"{a}hler, then $X$ is a projective and rationally connected manifold.}
\end{conjecture}

By employing the method of the above theorem 1.1 and the method of theorem 1.7 in \cite{NZ01} which use the form version of the Bochner identity, we prove the following theorem, which this partially confirms the above conjecture.

\begin{theorem}  Let $X$ be a compact complex manifold of complex demension $n>2$. Suppose $X$ has a Hermitian metric $h$ with positive real bisectional curvature, then the Hodge number $h^{1,0}=h^{2,0}=h^{n-1,0}=h^{n,0}=0$. In particular, if in addition $X$ is K\"{a}hler, then X is a projective. Also, it is rationallly connected manifold when $n=3$.
\end{theorem}

\begin{remark} In an attempt to generalize Wu-Yau's Theorem (\cite{WY}) to the Hermitian case, Yang and Zheng \cite{YZ} introduced the concept of {\em real bisectional curvature} for Hermitian manifolds. When the metric is K\"ahler, this curvature is the same as the holomorphic sectional curvature $H$, while when the metric is not K\"ahler, the curvature condition is slightly stronger than $H$ at least algebraically. This condition also appeared in the recent work by Lee and streets \cite{LS} where it is referred to as ``positive (resp. negative) curvature operator".
\end{remark}

What about compact K\"{a}hler manifolds with nonnegative holomorphic sectional curvature? Recently, S. Matsumura established some more general structure theorems for compact K\"{a}hler manifolds with semi-positive holomorphic sectional curvature (see \cite{Ma1},\cite{Ma2},\cite{Ma3}). In particular, he proved that any compact K\"{a}hler surface with semi-positive holomorphic sectional curvature is rationally connected, or a complex torus, or a ruled surface over an elliptic curve.
We might as well discuss the quasi-positive holomorphic sectional curvature, namely, nonnegative everywhere and positive at some point. In their paper \cite{HW}, Heier-Wong considered Yau's conjecture for projective varieties under assumption that the holomorphic sectional curvature is quasi-positive.

For complex dimension $n=2$, combine theorem 1.2 in \cite{Y2}, theorem 1.5 and remark 1.6 in \cite{HW}, we can easily get the following result.
\begin{theorem}
If compact complex manifold of dimension $n=2$ has a K\"{a}hler metric with quasi-positive holomorphic sectional curvature, then it is a projective and rationally connected manifold.
\end{theorem}

For high-dimensional situations, we can propose the following conjecture (see also \cite{Ma2}, Problem 4.3):
\begin{conjecture}{\em Let $X$ be a compact K\"{a}hler manifold with dimension $n\geq 3$. If $X$ has a
K\"{a}hler metric with quasi-positive holomorphic sectional curvature, then $X$ is a
projective and rationally connected manifold.}
\end{conjecture}

Similar to X. Yang's work on Yau'conjecture, L. Ni and F. Zheng (\cite{NZ01},\cite{NWZ}) have studied some new meaningful curvature, namely, orthogonal Ricci curvature. They proved that compact K\"{a}hler manifold with positive orthogonal Ricci curvature, then it must be projective. On the other hand, L. Ni and F. Zheng \cite{NZ02} also discussed the 2nd scalar curvature, which is the average of holomorphic sectional curvature over 2-dimensional subspaces of the tangent space. They proved that any compact K\"{a}hler manifold with positive 2nd scalar curvature must be projective. This is generalization of above mentioned result of X. Yang.

Because the method of X. Yang crucially relied on the geometric properties of RC-positivity, he also continued to delve into other properties of RC-positivity (\cite{Y3},\cite{Y4},\cite{Y5}).

For the negative holomorphic sectional curvature case, in the
recent breakthrough paper \cite{WY} of Wu and Yau, they proved that any projective
K\"{a}hler manifold with negative holomorphic sectional curvature must
have ample canonical line bundle. This result was obtained by Heier et. al.
earlier under the additional assumption of the Abundance Conjecture (\cite{HLW}).
In \cite{TY}, Tosatti and X.Yang proved that any compact K\"{a}hler manifold
with nonpositive holomorphic sectional curvature must have nef canonical
line bundle, and with that in hand, they were able to drop the projectivity
assumption in the aforementioned Wu-Yau Theorem. More recently, Diverio
and Trapani \cite{DT} further generalized the result by assuming that the
holomorphic sectional curvature is only quasi-negative, namely, nonpositive
everywhere and negative somewhere in the manifold. In \cite{WY1}, Wu and Yau
give a direct proof of the statement that any compact K\"{a}hler manifold with
quasi-negative holomorphic sectional curvature must have ample canonical
line bundle. We refer to \cite{HLW,WY,TY,DT,WY1,YZ,Nom,Tang} for more details.
\vspace{0.5cm}

\noindent\textbf{Acknowledgement.} The author is grateful to Professor Fangyang Zheng for constant encouragement and support. He wishes to express
his gratitude to Professor Shin-ichi Matsumura for useful discussions on \cite{Ma1},\cite{Ma2},\cite{Ma3}.
\vspace{0.5cm}

\section{Preliminaries}
In this section, we firstly summarize some basic properties about the RC-positivity which will be used in our proofs. In order to give geometric interpretations on rational connectedness, X. Yang \cite{Y1} introduced the following concept for Hermitian vector bundles:
\begin{definition}
Let $X$ be a complex manifold and $(E,h)$ be a Hermitian holomorphic vector bundle. Suppose $R^{(E,h)}\in \Gamma(X,\Lambda^{1,1}T_{X}^{\ast}\otimes End(E))$ is its Chern curvature tensor. $(E,h)$ is called $RC$-$positive$ (resp. $RC$-$negative$) if  any local nonzero section $a\in \Gamma(X,E)$,
there exists a local section $\nu\in \Gamma(X,T_{X})$, such that
\begin{align}
   R^{(E,h)}(\nu,\overline{\nu},a,\overline{a})>0.    (resp. <0)\nonumber
\end{align}
\end{definition}
For a line bundle $(L,h)$, it is RC-positive if and only if its Ricci curvature
has at least one positive eigenvalue at each point of $X$. This terminology
has many nice properties. For instances, quotient bundles of RC-positive
bundles are also RC-positive; subbundles of RC-negative bundles are still
RC-negative; the holomorphic tangent bundles of Fano manifolds can admit RC-positive metrics.
We refer to \cite{Y1} for more details.

In order to prove the Yau's conjecture, X. Yang first proved a main theorem.
\begin{theorem} [{\cite{Y1}, Theorem 1.4}] Let $X$ be a compact K\"{a}hler manifold of complex dimension $n$. Suppose that for every $1\leq p\leq n$,
there exists a smooth Hermitian metric $h_{p}$ on the vector bundle $\Lambda^{p}T_{X}$ such that $(\Lambda^{p}T_{X},h_{p})$ is RC-positive, then $X$ is
projective and rationally connected.
\end{theorem}

Another useful result is the following theorem:
\begin{theorem} [{\cite{Y1}, Corollary 6.6}] Let $X$ be a compact complex manifold. Suppose $X$ has a Hermitian metric with positive holomorphic
sectional curvature, then

(1) $T_{X}$ is RC-positive.

(2) $K_{X}^{-1}=\det T_{X}$ is RC-positive.
\end{theorem}

Let us recall the  concept of {\em real bisectional curvature} introduced in \cite{YZ}. Let $(X^{n},g)$ be a Hermitian manifold. Denote by $R$ the curvature tensor of the Chern connection. For $p\in X$, let $e=\{e_{1},\cdot\cdot\cdot,e_{n}\}$ be a unitary tangent frame at $p$, and let $a=\{a_{1},\cdot\cdot\cdot,a_{n}\}$ be non-negative constants with $|a|^{2}=a_{1}^{2}+\cdot\cdot\cdot+a_{n}^{2}>0$. Define the $real$ $bisectional$ $curvature$ of g by
\begin{align}
B_{g}(e,a)=\frac{1}{|a|^{2}}\sum_{i,j=1}^{n}R_{i\overline{i}j\overline{j}}a_{i}a_{j}.
\end{align}
We will say that a Hermitian manifold $(X^{n},g)$ has {\em positive real bisectional curvature}, denoted by $B_{g}>0$, if for any $p\in X$ and any unitary frame $e$ at $p$, any nonnegative constant $a=\{a_{1},\cdot\cdot\cdot,a_{n}\}$, it holds that $B_{g}(e,a)>0$.

Recall that the holomorphic sectional curvature in the direction $\upsilon$ is defined by $H(\upsilon)=R_{\upsilon\overline{\upsilon}\upsilon\overline{\upsilon}}/|\upsilon|^{4}$. If we take $e$ so that $e_{1}$ is parallel to $\upsilon$, and take $a_{1}=1$, $a_{2}=\cdot\cdot\cdot=a_{n}=0$, then $B$ becomes $H(\upsilon)$. So $B>0$ ($\geq0,<0,or \leq0$) implies $H>0$ ($\geq0,<0,or \leq0$).
For a more detailed discussion of this, we refer the readers to \cite{YZ}.

\vspace{0.5cm}

\section{Proof of Theorems 1.4}\label{tk,3}
\vspace{0.3cm}
\noindent {\bf Proof of Theorem of 1.4:} Since X has positive real bisectional curvature, then X has positive holomorphic sectional curvature, thus the canonical line bundle cannot admit any non-trivial global holomorphic section. In fact, its Kodaira dimension must be $-\infty$ which is proved by X.Yang in \cite{Y2}. So we just need prove $h^{1,0}=h^{2,0}=h^{n-1,0}=0$.

Let $s$ be a global holomorphic $p$-form on X, namely, $s\in H^{0}(X,\wedge^{p}\Omega)\cong H_{\overline{\partial}}^{p,0}(X)$. The Bochner identity gives
\begin{align}
\partial\overline{\partial}|s|^{2}=\langle \nabla s,\nabla s \rangle-\widetilde{R}(\cdot,\cdot,s,\overline{s})\nonumber
\end{align}
where $\widetilde{R}$ is the curvature of the Hermitian bundle $\Lambda^{p}\Omega$, and $\Omega=TX^{\ast}$ is the holomorphic cotangent bundle of X.

If $s$ is not identically zero, then $|s|^{2}$ will attain its nonzero maximum at some point $q$. Hence, at point $q$, we have
\begin{align}
\widetilde{R}(\upsilon,\upsilon,s,\overline{s})\geq 0\nonumber
\end{align}
For any type (1,0) tangent vector $\upsilon\in T_{q}X$. we want to show that this will contradict the assumption when $p=1$, $p=2$ or $p=n-1$.

When $p=1$, in a small neighborhood of $q$, we can write $s=\lambda \varphi_{1}$, where $\lambda\neq 0$ is a function and $\{\varphi_{1},\varphi_{2},\cdot\cdot\cdot,\varphi_{n}\}$ are local $(1,0)$-forms forming a coframe dual to a local tangent frame $\{e_{1},\cdot\cdot\cdot,e_{n}\}$, which is unitary at $q$. Since
\begin{align}
\widetilde{R}(\upsilon,\upsilon,s,\overline{s})=-|\lambda|^{2}R_{\upsilon\overline{\upsilon}1\overline{1}}\geq 0\nonumber
\end{align}
For any tangent direction $\upsilon$, where $R$ is the curvature tensor of $X$. If we take $\upsilon=e_{1}$, we have $R_{1\overline{1}1\overline{1}}\leq0$ which contradicts positive real bisectional curvature.

The $p=n-1$ case. In a small neighborhood of $q$, we can write $s=\lambda\varphi_{2}\wedge\cdot\cdot\cdot\wedge\varphi_{n}$, where $\lambda\neq 0$ is a function and $\{\varphi_{1},\varphi_{2},\cdot\cdot\cdot,\varphi_{n}\}$ are local $(1,0)$-forms forming a coframe dual to a local tangent frame $\{e_{1},\cdot\cdot\cdot,e_{n}\}$, which is unitary at $q$. Hence
\begin{align}
\widetilde{R}(\upsilon,\upsilon,s,\overline{s})=-|\lambda|^{2}\sum_{2}^{n}R_{\upsilon\overline{\upsilon}i\overline{i}}\geq 0\nonumber
\end{align}
for any $\upsilon\in T_{q}X$. If we take $\upsilon=e_{2},\cdot\cdot\cdot,e_{n}$, by adding them up, we get
\begin{align}
\sum_{i,j=2}^{n}R_{j\overline{j}i\overline{i}}\leq 0\nonumber
\end{align}
a contradiction to our assumption that real bisectional curvature is positive.

Now consider the $p=2$ case. Suppose that $s$ is non-trivial global holomorphic 2-form on X. Let $r\geq 1$ be the largest positive integer such that $s^{r}$ is not identically zero. Since we already have $h^{n,0}=h^{n-1,0}=0$, thus $2r\leq n-2$.

We will apply the Bochner formula to the $2r$-form $\beta=s^{r}$. Let $q$ be a maximum point of $|\beta|^{2}$. At point $q$, we can write $s=\sum_{i,j}a_{ij}\varphi_{i}\wedge \varphi_{j}$ under any unitary coframe $\{\varphi_{j}\}$ which is dual to a local unitary tangent frame $\{e_{j}\}$. The $n\times n$ is skew-symmetric. As is well-known (cf. \cite{Hua}), there exists unitary matrix $U$ such that ${}^{t}UAU$ is in the block diagonal form where each non-zero diagonal block is a constant multiple of $E$, where
\begin{align}
E=\left(
    \begin{array}{cc}
      0 & 1 \\
      -1 & 0 \\
    \end{array}
  \right)\nonumber
\end{align}
In other words, we can choose a unitary coframe $\varphi$ at point $q$ such that
\begin{align}
s=\lambda_{1}\varphi_{1}\wedge\varphi_{2}+\lambda_{2}\varphi_{3}\wedge\varphi_{4}+\cdot\cdot\cdot+\lambda_{k}\varphi_{2k-1}\wedge\varphi_{2k},\nonumber
\end{align}
where $k$ is a positive integer and each $\lambda_{i}\neq0$. Clearly, $k\leq r$ since $s^{k}\neq 0$ at $q$. If $k<r$, then $\beta=s^{r}=0$ at $q$, which is a maximum for $|\beta|^{2}$, implying $\beta=0$, a contradiction. Hence, we can get $k=r$. Thus $\beta=\lambda \varphi_{1}\wedge\varphi_{2}\cdot\cdot\cdot\wedge\varphi_{2r}$, where $\lambda=\lambda_{1}\cdot\cdot\cdot\lambda_{r}\neq 0$. By the Bochner formula, we get that
\begin{align}
\sum_{i=1}^{2r}R_{\upsilon\overline{\upsilon}i\overline{i}}\leq 0.\nonumber
\end{align}
If we take $\upsilon=e_{1},\cdot\cdot\cdot,e_{2r}$, by adding them up, we get
\begin{align}
\sum_{i,j=1}^{2r}R_{j\overline{j}i\overline{i}}\leq 0\nonumber
\end{align}
which contradicts the assumption that real bisectional curvature is positive.

If we assume that $X$ has a K\"{a}hler metric, then we have $h^{2,0}=h^{0,2}=0$. Hence, by Kodaira theorem (\cite{Ko}, Theorem 1), see also Proposition 3.3.2 and Corollary 5.3.3 \cite{Huy}, the $X$ is projective.

Finally, we will prove that it is rationallly connected manifold when $n=3$. By Theorem 2.2, it is enough to show that for every $1\leq p\leq 3$, $\Lambda^{p}TX$ is $RC$-positive. Since positive holomorphic sectional curvature implys that $TX$ and $K_{X}^{-1}=\det TX$ are $RC$-positive which is from Theorem 2.3, we just need prove $\Lambda^{2}TX$ is $RC$-positive or $\Lambda^{2}\Omega$ is $RC$-negative.

Suppose $(\Lambda^{2}\Omega, \Lambda^{2}h^{\ast})$ is not $RC$-negative. By Definition of $RC$-negative, there exist a point $q\in X$ and a nonzero local section $a\in \Gamma(X,\Lambda^{2}\Omega)$ such that
\begin{align}
\widetilde{R}(\upsilon,\overline{\upsilon},a,\overline{a})\geq 0 \nonumber
\end{align}
For any type (1,0) tangent vector $\upsilon\in T_{q}X$. In a small neighborhood of $q$, we can write $a=\lambda\varphi_{2}\wedge\varphi_{3}$, where $\lambda\neq 0$ is a function and $\{\varphi_{1},\varphi_{2},\varphi_{3}\}$ are local $(1,0)$-forms forming a coframe dual to a local tangent frame $\{e_{1},e_{2},e_{3}\}$, which is unitary at $q$. Hence
\begin{align}
\widetilde{R}(\upsilon,\upsilon,a,\overline{a})=-|\lambda|^{2}\sum_{2}^{3}R_{\upsilon\overline{\upsilon}i\overline{i}}\geq 0\nonumber
\end{align}
for any $\upsilon\in T_{q}X$. If we take $\upsilon=e_{2},e_{3}$, by adding them up, we get
\begin{align}
\sum_{i,j=2}^{3}R_{j\overline{j}i\overline{i}}\leq 0\nonumber
\end{align}
a contradiction to our assumption that real bisectional curvature is positive. So, we have completed the proof of theorem 1.4. \qed

\end{document}